\documentclass{amsart}

\usepackage{amsmath,amscd, amssymb}
\usepackage{colordvi}

\def\Var{\mathrm{Var}}
\def\Lie{\mathrm{Lie}}
\def\RVar{\mathrm{RVar}}
\def\RLie{\mathrm{RLie}}
\def\RAVar{\mathrm{RAVar}}
\def\RALie{\mathrm{RALie}}
\def\RBLie{\mathrm{RBLie}}
\def\preCom{\mathrm{preCom}}

\def\CYRSH{\font\cyr=cmcyr10\cyr \char'173}

\let\<\langle
\let\>\rangle

\def\LS{\mathop{\fam0 LS}\nolimits}
\def\RLS{\mathop{\fam0 RLS}\nolimits}
\def\RALS{\mathop{\fam0 RALS}\nolimits}
\def\Jac{\mathop{\fam0 Jac}\nolimits}
\def\Ker{\mathop{\fam0 Ker}\nolimits}
\def\End{\mathop{\fam0 End}\nolimits}
\def\gr{\mathop{\fam0 gr}\nolimits}

\newtheorem{lemma}{Lemma}

\newtheorem{theorem}{Theorem}
\newtheorem{corollary}{Corollary}
\newtheorem{remark}{Remark}
\newtheorem{proposition}{Proposition}
\newtheorem{example}{Example}

\title[Rota---Baxter Lie algebras]{Gr\"obner---Shirshov basis of the universal enveloping Rota---Baxter algebra\\ of a Lie algebra}
\thanks{Supported by Russian Science Foundation (project 14-21-00065)}
\author{Vsevolod Gubarev, Pavel Kolesnikov}
\address{Sobolev Institute of Mathematics, Novosibirsk, Russia}

\begin{document}


\maketitle
 
\section{Introduction}

A Rota---Baxter algebra is a linear space $A$ over a field $\Bbbk $
equipped with bilinear product $(a,b)\mapsto ab$, $a,b\in A$, 
and with a linear map $R: A\to A$ such that 
\begin{equation}\label{eq:RB_def}
 R(a)R(b) = R(R(a)b) + R(aR(b)) + \lambda R(ab),
\end{equation}
where $\lambda $ is a constant from $\Bbbk $.
A linear operator $R$ satisfying \eqref{eq:RB_def} is called a Rota---Baxter operator of weight $\lambda $. 

This notion initially appeared in analysis \cite{Baxter60}, and then in combinatorics \cite{Rota69} and quantum field theory 
\cite{ConnKreimer2000}. We refer the reader to the book \cite{Guo_RBA} and references therein for more details.
There is a number of studies on associative and commutative Rota---Baxter algebras. Let us mention those that 
are close to the topic of this paper.

A linear basis of the free associative Rota---Baxter algebra was found in 
\cite{FardGuo2008}, where it was also shown that the universal enveloping 
Rota---Baxter algebra of a free dendriform (or tridendriform, for nonzero weight) algebra is free. 
A simpler proof 
of the same fact follows from \cite{GubKol2013}. Another method for finding this basis 
was applied in \cite{Bokut2010}, in a more modern form this approach was exposed in \cite{Guo}.

The class of Rota---Baxter Lie algebras is of special interest since it is closely related 
with pre-Lie (left/right-symmetric) algebras. Namely, if $L$ is a Lie algebra with
a product $[\cdot , \cdot]$ equipped with a Rota---Baxter operator
$R$ then the same space $L$ with new operation $ab = [R(a),b]$, $a,b\in L$, 
is a pre-Lie algebra. 

Moreover, there is a natural relation between Rota---Baxter operators and 
solutions of the classical Yang---Baxter equation (CYBE) \cite{Semenov83}. Namely, 
if $L$ is a Lie algebra equipped with 
an symmetric invariant bilinear form $\<\cdot ,\cdot \>$ (not necessarily non-degenerate) 
then there is a natural map
$L\to L^*$, $a\mapsto \<a,\cdot\>$, and thus we have a map
$\Phi : L\otimes L \to \End(L)$.
If $X\in L\otimes L$ is a skew-symmetric solution of CYBE
\[
 [X^{12},X^{13}] + [X^{12},X^{23}] + [X^{13},X^{23}]=0.
\]
then  $R=\Phi(X)$ is a Rota---Baxter operator on $L$. 

This paper is devoted to combinatorial structure of Lie algebras with a Rota---Baxter operator. 
The main problem we solve is an analogue of the PBW Theorem for universal enveloping 
Rota---Baxter Lie algebra $U_{RB}(L)$ of an arbitrary Lie algebra $L$. 
We prove that $U_{RB}(L)$ carries natural filtration such that the corresponding associated graded algebra 
$\gr U_{RB}(L)$ is isomorphic (as a Lie algebra) to the universal enveloping one in the class $\RALie $ of Lie algebras 
with linear operator $R$ satisfying the identity $[R(x),R(y)]=0$.
We also note that the same statement is true in the varieties As and Com of associative and commutative algebras.

The main tool of the proof is a version of the Composition-Diamond Lemma (CD-Lemma) 
for Lie algebras with an additional operator. 
A more general approach to this Lemma (for Lie algebras with an arbitrary set of additional operators) 
was developed in \cite{ChenQiu16}.
In the proof of CD-Lemma we use terminology of \cite{Guo} and some combinatorial results 
of \cite{Sh58_free_lie}. For more detailed exposition of the latter results, see \cite{BokChen13}.

\section{Algebras with additional operator}

Suppose $\Var $ is a variety of linear algebras. 
Denote by $\RVar $ the variety of $\Var$-algebras equipped with 
an additional linear operator $R$.
Denote the natural forgetful functor from $\RVar$ to $\Var $ by $\Theta_R $,
and let $U_R$ stands for its left adjoint functor from $\Var $ to $\RVar $.


Obviously, for the free $\Var$-algebra $\Var \<V\>$ generated by a linear space $V$
the universal $\RVar$-envelope $U_R(\Var\<V\>)$ is isomorphic to 
the free $\RVar$-algebra $\RVar\<V\>$.

Let us state the explicit construction of $U_R(A)$.
Given $A\in \Var $, denote by 
$\bar A$ the copy of the linear space $A$. Let $\rho : A\to \bar A$
stands for the isomorphism $a\to \bar a$.

Construct a series of algebras $\{A_n\}_{n\ge 0}$ by the following rule:
\[
\begin{gathered}
A_0=A  , \\
A_1 = A* \Var\<\bar A \>, \\
\hdots \\
A_n = A* \Var\< \bar A_{n-1}\> , \\
\hdots
 \end{gathered}
\]
where $*=*_\Var $ denotes the free product in the variety $\Var $.

As above, let $\bar A_n$ denotes a copy of the space $A_n$; denote the linear isomorphism 
$a\to \bar a$, $a\in A_n$, by $\rho_n$.

Construct a series of $\Var $-homomorphisms $\tau_n: A_{n}\to A_{n+1} $, $n\ge 0$, as follows.
Set $\tau_0$ be the canonical embedding of $A$ into the free product $A_1$ and proceed by induction:
\[
\begin{CD}
A_{n-1} @>\tau_{n-1}>> A_n \\
@VV\rho_{n-1} V  @VV\rho_n V \\
\bar A_{n-1} @>\bar \tau_{n-1}>> \bar A_n \\
@VV\subset V  @VV\subset V \\
\Var\<\bar A_{n-1}\> @>\tau_{n-1}^0>> \Var\<\bar A_n\> \\
@VV\subseteq V  @VV\subseteq V \\
A*\Var\<\bar A_{n-1}\> @>\Red{\tau_{n}}>> A*\Var\<\bar A_n\>,
\end{CD}
\]
where $\bar \tau_{n-1}=\rho_{n-1}^{-1}\circ \tau_{n-1}\circ \rho_n$
is a linear homomorphism, $\tau_{n-1}^0$ is the induced $\Var $-homomorphism 
of free algebras, and $\tau_n$ is a $\Var$-homomorphism that comes 
from the definition of free product.

\begin{lemma}\label{lem:SpectrumT}
All homomorphisms $\tau_n$ are injective.
\end{lemma}

\begin{proof}
Assume $\tau_{n-1}$ is injective. 
Consider
\[
\begin{CD}
 \bar A_{n-1} @>\bar\tau_{n-1}> \subseteq > \bar A_n @>\varphi >> \bar A_{n-1}\\
 @V\subseteq VV  @V\subseteq VV @ VVV\\
 \Var\<\bar A_{n-1} \> @>>> \Var \<\bar A_n\> @>>> \Var\<\bar A_{n-1}\>
 \end{CD}
\]
where $\varphi $ is a projection of the linear space $\bar A_{n}$ onto $\bar A_{n-1}$: $\varphi\tau_{n-1}^0=id_{\bar A_{n-1}}$
(every such linear map extends to a homomorphism of free algebras). 
Then the universal property of free product (uniqueness) implies 
the existence of 
\[
\begin{CD}
 A_n=A*\Var\<\bar A_{n-1}\> @>\tau_{n}>> A*\Var\<\bar A_{n}\> @>\Red{\psi} >> A*\Var\<\bar A_{n-1}\>=A_n,
 \end{CD}
\]
where 
$\psi \tau_n =id_{A_n}$, so $\tau_n$ is also injective.
\end{proof}

\begin{lemma}\label{lem:SpectrumR}
For every $\RVar$-algebra $B$ and for every homomorphism of $\Var $-algebras $\psi: A\to B$
there exists unique family $\{\psi_n\}_{n\ge 0}$ of $\Var$-homomorphisms $\psi_n: A_n\to B$ such that 
\[
 \rho_n\circ \psi_{n+1} = \psi_n\circ R
\]
and
\[
 \psi_0=\psi,\quad \psi_n = \tau_n\circ \psi_{n+1}.
\] 
\end{lemma}

\begin{proof}
Let us show existence and uniqueness by induction. 
Given $\psi_n: A_n\to B$, construct 
\[
\begin{CD}
 A_n @>\rho_n>> \bar A_n @>\subseteq >> \Var\<\bar A_n\> @>\subseteq >> A_{n+1}\rlap{${}=A*\Var\<\bar A_n\>$} \\
 @V\psi_n VV  @V\bar \psi_n VV @V\psi_n^0 VV @AA\subseteq A \\ 
 B @>R>> B @>\mathrm{id} >> B @<\psi << A
\end{CD}
\]
Here the rightmost vertical arrow 
is the canonical embedding of $A$ into the free product which coincides with 
$\tau_0\circ \dots \circ \tau_n$, $\bar\psi_n=\rho_n^{-1}\circ \psi_n\circ R$ is a linear map, 
$\psi_n^0$ is a homomorphism of $\Var $-algebras induced by $\bar \psi_n$.
The right-hand square in the diagram above induces $\Var $-homomorphism $\psi_{n+1}: A*A_n^0 \to B$.

Why $\psi_n = \psi_{n+1}\tau_n$? For $n=0$, it follows from the definition of $\psi_1$.
Assume $n>0$ and $\psi_{n-1}=\psi_n\tau_{n-1}$. 
Then for all $y\in \bar A_{n-1}$
\[
 \bar\tau_{n-1}(y) = \rho_n\tau_{n-1}\rho_{n-1}^{-1}(y)
\]
Since $\bar\psi_n\rho_n(z) = R\psi_n(z)$ for all $z\in A_n$, we have
\[
 \bar\psi_n\bar\tau_{n-1}(y) 
 = \bar\psi_n\rho_n\tau_{n-1}\rho_{n-1}^{-1}(y)
 =R\psi_n\tau_{n-1}\rho_{n-1}^{-1}(y)
 =R\psi_{n-1}\rho_{n-1}^{-1}(y)
 =\bar\psi_{n-1}(y).
\]
Therefore, the induced $\Var $-homomorphisms are related in the same way:
\[
 \psi^0_{n-1}=\psi_n^0\tau_{n-1}^0.
\]
Now, for all $x\in \Var\<\bar A_{n-1}\> \subseteq A_n$ we have 
\[
 \tau_n(x) = \tau_{n-1}^0(x).
\]
Hence,
\[
 \psi_{n+1}\tau_n(x)
  = \psi_n^0\tau_{n-1}^0(x) = \psi_{n-1}^0(x) = \psi_n(x)
\]
by definition of $\psi_{n+1}$. 
Since $\psi_{n}$ is uniquely determined by its action on $\bar A_{n-1}$ 
(uniqueness property of the universal map on free product), 
we have the required equality $\psi_{n+1}\tau_n = \psi_n$ on the entire~$A_n$. 
\end{proof}

The chain 
\[
 \begin{CD}
  A @>\tau_0 >> A_1 @>\tau_1 >> A_2 @>>> \dots  @>>> A_{n}  @>\tau_n >> A_{n+1} @>>> \dots  
 \end{CD}
\]
naturally defines direct system of $\Var$-algebras. 
Let
\[
 A_\infty = \lim\limits_{\to} A_n,
\]
\[
 \rho : A_\infty \to A_\infty ,\quad \rho = \lim\limits_{\to}\rho_n.
\]

\begin{theorem}\label{thm:REnvelope}
 The $\Var$-algebra $A_\infty $ with linear map $\rho $  
 is isomorphic to the universal $\RVar$ enveloping $U_R(A)$.
\end{theorem}

\begin{proof}
 The universal property of $(A_\infty, \rho )$ follows from Lemma~\ref{lem:SpectrumR}.
\end{proof}

Let us consider the particular case $\Var=\Lie$. Recall that if $Y$ is a well-ordered set of 
generators then the linear basis of $\Lie\<Y\>$ may be constructed in the following way \cite{Sh58_free_lie}. 
A word $u\in Y^*$ is called an (associative) {\em Lyndon---Shirshov word} (LS-word) if either $u\in Y$ or for every 
presentation $u=vw$, $v,w\in Y^*$, we have $u>wv$ lexicographically. 
Denote the set of all such words by $\LS(Y)$.
For every $u\in \LS(Y)$ there exists {\em standard} bracketing $[u]$ 
such that $[u]=([v][w])$, where $w$ is the longest proper LS-suffix of $u$ 
(then $v$ is also an LS-word, $[v]$ and $[w]$ are standard bracketings on these shorter words).
The set $\{[u] \mid u\in\LS(Y) \}$ is a linear basis of $\Lie \<Y\>$.

It is not hard to construct the linear basis of free R$\Lie$-algebra $\RLie\<X\>$
for a given well-ordered set $X$ of generators.
Let
 $\RLS_0(X)=\{[u] \mid u\in \LS(X) \}$ be the basis of $\Lie \<X\>$ equipped with leg-lex ordering: 
\[
 [u]<[v] \iff u<_{deglex} v
\]
Assume the set $RLS_n(X)$ is already constructed and equipped with a well order.

Consider the alphabet $U_n = X\cup \{R([u]) \mid [u]\in \RLS_n(X) \}$ with the following
order: 
$x<R([u])$ for all $x\in X$, $[u]\in \RLS_n(X)$; $R([u])<R([v]) \iff [u]<[v]$, $[u],[v]\in \RLS_n(X)$.
Then 
\[
 \RLS_{n+1}(X):=\{[w] \mid w\in \LS(U_n) \}
\]
equipped with deglex order.

Obviously, $\RLS_n(X)\subset \RLS_{n+1}(X)$ for all $n\ge 0$.
%
%

\begin{corollary}\label{cor:RLie}
The set 
\[
 \RLS(X) = \bigcup\limits_{n\ge 0} \RLS_n(X)
\]
is a linear basis of $\RLie\<X\>$.
\end{corollary}

\begin{proof}
Let $L=\Lie\<X\>$, $U_R(L)\simeq \RLie \<X\> \overset{\theta}{\to} L_\infty $.
Consider the images of RLS-words as elements of $\RLie\<X\>$ under the isomorphism $\theta $ iduced by $x\mapsto x$, $x\in X$.
By definition, $\theta(\RLS_0(X))$ is the basis of $L_0=L$. 
Assume $\theta(\RLS_k(X))$ is a basis of $L_k$ for  all $k\le n$, and the embedding 
$\RLS_{k-1}(X)\subset \RLS_k(X)$ is compatible with $\tau_{k-1}:L_{k-1}\to L_k$, $k=1,\dots, n$.
Then $\theta(R(\RLS_n(X)))=\rho_n(\theta(\RLS_n(X)))$ is the set of free generators for 
$\Lie\<\bar L_n\>$. Moreover, $\theta $ is compatible with $\tau_{n+1}$.

Recall that a linear basis of a free product of two free Lie algebras 
is the free Lie algebra generated by disjoint union of the generating sets.
In our case, one of these sets is $X$, other is $\theta(R(\RLS_n(X)))$.
Therefore, $\theta(\RLS_{n+1}(X))$ is the linear basis of $L_{n+1}$.
\end{proof}

In particular, 
$\RLie \<X\> $ as a Lie algebra is isomorphic to $\Lie \<U\>$, 
where $U=\bigcup\limits_{n\ge 0} U_n$. 
Therefore, $\RLie \<X\>$ has a natural ascending filtration 
\begin{equation}\label{eq:Deg-Filtration}
\RLie^{(n)}\<X\> =\{ f\in \RLie\<X\> \mid \deg f\le n \},
\end{equation}
where $\deg f$ is the degree of $f\in \Lie\<U\>$ relative to the alphabet~$U$.

Note that $U$ may not be a well-ordered set, e.g., 
$xy>R(xy)>R^2(xy)> \dots $ for $x>y$. 
However, for every $n\ge 0$ the subset $U_n$ is obviously well-ordered.

For an RLS-word $[u]$, denote by $\deg_R u$ ($R$-degree) the total number of operators $R$ appearing in $u$.
For $f\in \RLie \<X\>$, set $\deg_R f$ to be the maximal $R$-degree among all its monomials.
Note that for every $n\ge 0$ there exists $N$ such that  
$\{[u]\in \RLS(X) \mid \deg_R u\le n \}\subseteq \LS(U_N)$.

\begin{remark}\label{rem:RAVar}
Denote by $\RAVar$ the subvariety of $\RVar$ defined by identity $R(x)R(y)=0$
(image of $R$ is abelian).  The following construction is completely similar 
to the one stated above. 
\end{remark}

For $A\in \Var$, let $A_0=A$ and $A_{n+1} = A*A_n^0$, $n\ge 0$, where $A_n^0$ stands 
for the same space as $A_n$ considered as an algebra with trivial operations.
Then the universal enveloping $\RAVar$-algebra $U_{RA}(A)$ is isomorphic 
to $\lim\limits_{\to} A_n$.

\begin{corollary}
 The free $\RALie$-algebra $\RALie \<X\>$ is isomorphic as a Lie algebra to 
 the partially commutative Lie algebra
 $\Lie \<U \mid uv=0,\ u,v\in U\setminus X \> $.
\end{corollary}

Let us denote by $\RALS_n(X)$ the set of all $[w]\in \RLS_n(X)$ such that $w$ do not contain subwords
of the form $R([u])R([v])$, $u,v\in \RLS_{n-1}(X)$, $[u]>[v]$. 
It is easy to see \cite{Sh62_one_hypoth} that $\RALS_n(X)$ is the linear basis of the Lie algebra $L_n$ 
constructed from $L_0=\Lie\<X\>$ as above. Therefore, 
\[
\RALS(X) = \bigcup\limits_{n\ge 0} \RALS_n(X) 
\]
is the linear basis of $U_{RA}(\Lie\<X\>)\simeq \RALie \<X\>$.

\section{CD-lemma for RLie algebras}

Let us call elements of $\RLie\<X\>$ by RLie-polynomials, and let 
$\bar f \in \RLS(X) $ stand for the leading word (principle monomial) of an RLie-polynomial~$f$.

Let us recall an important statement which plays an important role in the combinatorial theory 
of Lie algebras.

\begin{lemma}[{Shirshov bracketing, \cite[Lemma~4]{Sh58_free_lie}}]
Let $U$ be an ordered set, and $w,u\in \LS(U)$.
Suppose $u$ is a subword of $w$, i.e., $w=aub$, where $a$ and $b$ are some words in $U$
(either of them may be empty). Denote by $w_{u\leftarrow *} = a{*}b$, a word in the alphabet 
$U\dot\cup \{*\}$ obtained from $w$ by replacing this occurence of $u$ by a new symbol~$*$.
Then there exists unique bracketing on $w_{u\leftarrow *}$, denoted by $\{w_{u\leftarrow *}\}$,
such that 
\[
 \{a[u]b\} = [w] +\sum\limits_i \alpha_i [w_i],\quad \alpha_i\in \Bbbk,\ w_i\in \LS(U), \ [w_i]<[w].
\]
\end{lemma}

Uniqueness of the Shirshov bracketing implies the following property:
let $w,u,z\in \LS(U)$, $u$ is a subword of $z$, and $z$ is a subword of $w$.
Consider the words 
 $w_{z\leftarrow *} = a{*}b$, $w_{u\leftarrow *}$, and $z_{u\leftarrow *}$ 
 with the corresponding Shirshov bracketings $\{\dots \}$. Then 
\[
 \{a \{ z_{u\leftarrow *} \} b\} = \{ w_{u\leftarrow *} \}.
\]

Suppose $S$ is a set of monic RLie polynomials. Construct $\hat S$ as follows.
For every $f\in S$, $\bar f=[u]$, consider the associative word $u\in \LS(U)$ and consider all $[w]\in \RLS(X)$ such that 
the corresponding $w\in \LS(U)$ contain $u$ as a subword: $w=aub$.
Let $\{w_{u\leftarrow *}\} = \{a{*}b\}$ be the Shirshov bracketing.
Denote by $\hat S_0$ the collection of all RLie polynomials $\{w_{u\leftarrow f}\} = \{afb\}$
corresponding to all possible occurences of $u$, $[u]=\bar f$, $f\in S$, in all RLS-words $[w]$.
Then $\overline{w_{u\leftarrow f}}  =[w]$ and $w_{u\leftarrow f} $ belongs to the ideal of the Lie algebra 
$\Lie\<U\>$ generated by $S$.
All these polynomials are monic, and $S\subset \hat S_0$.

Proceed by induction: given $\hat S_n = \Sigma$, define $\hat S_{n+1} = \Sigma \cup \widehat{R(\Sigma )}_0\supset \hat S_n$, 
and
\[
 \hat S = \bigcup\limits_{n\ge 0} \hat S_n
\]

\begin{lemma}\label{lem:R-Ideal_sum}
 An RLie polynomial $f$ belongs to the ideal $I_R(S)$ generated by $S$ in $\RLie \<X\>$ if and only if 
$f  = \sum_i \alpha_i h_i$, $h_i\in \hat S$, $\alpha_i\in \Bbbk $.
\end{lemma}

\begin{proof}
 The ideal $I_R(S)$ in $\RLie\<X\>$ is the minimal $R$-invariant ideal in the Lie algebra $\Lie\<U\>$ which contains $S$. 
 By the construction, $I_R(S)\subseteq \hat S$. 
 
 Conversely, it follows from \cite[Lemma~3]{Sh62_alg_Lie} that an ideal $I(\Sigma )$ generated by a set $\Sigma $ in $\Lie \<U\>$
 coincides with the linear span of $\hat \Sigma_0$. Hence, the linear span of $\hat S$ is an ideal in $\Lie\<U\>$. Obviously, 
 this ideal is $R$-invariant, so $I_R(S)\subseteq \Bbbk \hat S$.
\end{proof}

Recall that a {\em rewriting system}  is an oriented graph $\mathcal G=(V,E)$ which has no infinite oriented paths.
A vertex $v\in V$ is called {\em terminal} if there are no edges of the form $v\to w$ in $E$.

Define an oriented graph $\mathcal G_R(X,S)$ on the set of vertices
$\RLie \<X\>$ based on a set of monic RLie-polynomials $S$, assuming that two RLie-polynomials $f$ and $g$ are connected by an edge 
$f\to g$ if and only if 
$f = f_0 + \alpha [u] + f_1$ (all monomials of $f_0$ are larger than $[u]$ and $[u]>\bar f_1$, $\alpha\in \Bbbk $, $\alpha \ne 0$)
such that $[u]= \bar h$ for some $h \in \hat S$, and $g = f - \alpha h$.
Every edge obviously corresponds to unique $h\in \hat S$, and therefore has a well-defined {\em level} which is the minimal
  $n$ such that $h\in \hat S_n$.

From now on, assume the following additional condition on $S$: 
$\deg_R \bar s \ge \deg_R s$
for every $s\in S$, i.e., 
the number of operators $R$ in the leading word $\bar s$
is greater or equal to $R$-degrees of all other monomials in~$s$.
Obviously, the same relation holds for $h\in \hat S$.
In this case, $\mathcal G_R(X,S)$ is a rewriting system since for every vertex $f\in \RLie \<X\>$ its cone (set of all vertices $g$ such that 
there exists an oriented path $f\to \dots \to g$) belongs to $\Bbbk\RLS_n(X)$ for some $n$, and 
$U_n^*$ is well ordered.
Terminal vertices of this rewriting system are also called {\em $S$-reduced} RLie polynomials.

Let us denote by $f\sim _d g$ the fact that $f,g\in \RLie\<X\>$ are connected by a non-oriented path of length $d\ge 1$.
Notation $f\sim g$ means that there exists $d\ge 1$ such that $f\sim_d g$.

The following two lemmas are almost obvious but we still state their proofs for 
readers' convenience.

\begin{lemma}\label{lem1}
Let $V$ be a subspace of $\RLie\<X\>$, and let 
$\mathcal G(V)$ stand for the subgraph of $\mathcal G_R(X,S)$ with vertices $V$.
Then for every $f,g,h\in V$ 
 \[
 f\sim g \text{ in $\mathcal G(V)$} \iff f+h\sim g+h \text{ in $\mathcal G(V)$}. 
 \]
\end{lemma}

\begin{proof}
 It is enough to show $(\Rightarrow)$. Suppose $f\sim_d g$ and proceed by induction in $d$.
 In fact, we only need $d=1$ since the induction step is obvious.
 Assume $f\to g$,
 $f = f_0 + \alpha [u] + f_1$, $[u]= \bar s$, $s \in \hat S$, 
 $g = f - \alpha s$ as in the definition of $\mathcal G_R(X,S)$.
 In particular, $s\in V$
 Apply the same principle to write down decompositions of $h = h_0 +\beta [u] + h_1$ (for some $\beta \in \Bbbk $)
 and $g = f-\alpha s = g_0+g_1$.
Then 
$f+h = f_0+h_0 + (\alpha +\beta )[u] + f_1+h_1$, 
$g+h = g_0+h_0 +\beta [u] + h_1 + g_1$.
If $\alpha +\beta\ne 0$ and $\beta = 0$ then $f+h\to g+h$.
If $\alpha +\beta\ne 0$ and $\beta \ne 0$ then 
\[
 f+h\to f + h - (\alpha+\beta)s = g + h-\beta s \leftarrow g+h,
\]
so $f+h\sim_2 g+h$. 
Finally, if $\alpha+\beta=0$ then 
\[
 f+h = f_0+h_0 +f_1+h_1 = g +\alpha s +h = g+h-\beta s \leftarrow g+h.
\]
\end{proof}

\begin{lemma}\label{lem2}
In the notations of Lemma~\ref{lem1}, the following statement holds:
for every $f,g\in V$
\[
f-g = \sum\limits_{i} \alpha_i s_i, \quad \alpha_i\in \Bbbk ,\ s_i \in \hat S\cap V,
\]
if and only if $f\sim g$ in $\mathcal G(V)$.
\end{lemma}

\begin{proof}
($\Leftarrow $) It follows from the definition of edges in $\mathcal G_R(X,S)$. 

($\Rightarrow $)
Assume $f-g=\alpha_1 s_1+ \dots + \alpha_n s_n$, $s_i\in \hat S\cap V$.
If $n=1$ then we simply have $f-g\to 0$, so $f\sim g$ by Lemma \ref{lem1}.
If $n>1$ then $f-(g+\alpha_1 s_1)\sim 0$ by induction, so $f-g\sim \alpha_1s_1\to 0$ by Lemma~\ref{lem1}.
\end{proof}

By Lemma~\ref{lem:R-Ideal_sum}, the ideal $I_R(S)$ coincides with the linear span of $\hat S$.
Therefore, connected components (in the non-oriented sense) of $\mathcal G_R(X,S)$ are exactly the elements 
of the quotient algebra $\RLie\<X\>/I_R(S)$.

We will mainly use the following subspaces of $\RLie \<X\>$: 
\[
\begin{gathered}
 V_n = \Bbbk \{[u]\in \RLS(X) \mid \deg_R u \le n\},\quad n\ge 0, \\
 V^{[w]} = \Bbbk \{[u]\in \RLS(X) \mid [u]\le [w]\}, \quad [w]\in \RLS(X), \\
 V_n^{[w]} = V_n\cap V^{[w]}.
\end{gathered}
\]
Note that 
\[
V_n=\bigcup\limits_{[w]\in \RLS(X)\cap V_n} V_n^{[w]},
\]
and $\RLS(X)\cap V_n $ is a well-ordered subset of $\RLS(X)$.

Recall that a rewriting system is called {\em confluent} if for every vertex $v$ there exists 
unique terminal vertex $t$ such that $v$ is connected 
with $t$ by an oriented path (i.e., $v\to \dots \to t$, or $v\rightsquigarrow t$). 
In particular, every non-oriented connected component of a confluent rewriting system 
contains unique terminal vertex. 

Therefore, if the rewriting system $\mathcal G_R(X,S)$ is confluent then there exists unique 
normal form of an element of $\RLie\<X\>/I_R(S)$ which may be found by straightforward 
walk on the graph. The following statement is a well-known 
criterion of confluentness.

\begin{theorem}[Diamond Lemma, \cite{Newm53}]\label{thm:Diamond}
 A rewriting system $\mathcal G=(V,E)$ is confluent if and only if 
 for every $v\in V$ and for every two edges $v\to w_1$, $v\to w_2$
 there exists a vertex $u\in V$ such that $w_1\rightsquigarrow u$ and $w_2 \rightsquigarrow u$. 
  \qed
\end{theorem}

It is easy to see that rewriting system $\mathcal G_R(X,S)$ is confluent if and only if 
so is each subsystem $\mathcal G(V_n)$, $n\ge 0$. 
The latter is confluent if and only if so is $\mathcal G(V_n^{[w]})$, 
$[w]\in \RLS(X)\cap V_n$.

\begin{proposition}\label{prop:Ambiguity1}
Let $S\subset \RLie \<X\> $ be a set of monic RLie-polynomials, $n\ge 0$.
Suppose the rewriting system 
$\mathcal G(V_n)\subset \mathcal G_R(X,S)$
has the following property:
for every RLS-word  $[w]\in V_n$ and for every pair of edges $[w]\to g_1$, $[w]\to g_2$ in $\mathcal G(V_n)$
we have 
\begin{equation}\label{eq:TrivialComp}
g_1-g_2 = \sum\limits_i \alpha_i h_i,\quad h_i\in \hat S\cap V_n,\ \bar h_i<[w].
\end{equation}
Then the system $\mathcal G(V_n)$ is confluent.
\end{proposition}

\begin{proof}
Let us check the Diamond Condition from Theorem \ref{thm:Diamond}
for rewriting system 
$\mathcal G_n^{[v]}=\mathcal G(V_n^{[v]})\subset \mathcal G_R(X,S)$, 
$[v]\in \RLS(X)\cap V_n$. 

Proceed by induction on $[v]$. Assume 
the rewriting system $\mathcal G_n^{[u]}$ is confluent for all $[u]\in V_n$, $[u]<[v]$, 
and consider an ambiguity in the graph $\mathcal G_n^{[v]}$, 
i.e., a pair of edges 
$f\to g_1$, $f\to g_2$.
Here $g_1 = f - \alpha h_1$, $g_2 = f - \beta h_2$, $h_i\in \hat S\cap V_n^{[v]}$.
There are three possible cases:

Case~1: $\bar h_1\ne \bar h_2$. \\
Then $f$ may be written in the form with ordered monomials 
$f = f_0 +\alpha [u_1] + f_{01} + \beta [u_2] + f_1$, $[u_i]=\bar h_i$.
Suppose $[u_1] - h_1 = \gamma [u_]2 + h$, where $h$ does not contain monomial $[u_2]$, 
$\gamma \in \Bbbk $. 
It is now easy to see that if 
$\gamma\alpha + \beta \ne 0$ then there exists edges $g_1\to g_1'\to g$, $g_2\to g$
in $\mathcal G_n^{[v]}$, where $g= f_0+f_{01}+f_1 +\alpha h +(\beta + \gamma\alpha )([u_2]-h_2)$. 
If $\gamma\alpha +\beta =0$ then there exist edges $g_1\to g$, $g_2\to g_2'\to g$ 
for the same~$g$. Therefore, in this case the Diamond Condition holds.

Case 2: $\bar h_1=\bar h_2<\bar f$.\\
Then $f  = f_0 +\alpha [u] + f_1$, $[u]=\bar h_i<\bar f\le [v]$.
Hence, there in an ambiguity 
$f'\to g_1'$, $f'\to g_2'$ in $\mathcal G_n^{[u]}$, where 
$f'=\alpha [u]+f_1$, $g_i= f_0 + g_i'$. 
By the inductive assumption, 
there exist two paths in $\mathcal G_n^{[u]}$:
$g_i'\to \dots \to g'$, $i=1,2$. Therefore, 
$g_i\to \dots \to f_0+g'$ in $\mathcal G_n^{[v]}$ since all monomials in $f_0$ 
are greater than $[u]$.

Case 3: $\bar h_1=\bar h_2=\bar f$.\\
Without loss of generality, assume $\bar f=[v]$.
Then $g_i = \alpha ([v]-h_i) + f_1$ and the difference
$g_1-g_2 = \alpha (h_2-h_1)$ coincides (up to scalar) with one that appears 
in the pair of edges $[v]\to [v]-h_i$, $i=1,2$. 
Therefore, the condition of the statement implies
$g_1$ and $g_2$
are connected by a non-oriented path in $\mathcal G_R^{[u]}(X,S)$
for some $[u]<[v]$. The last rewriting system is assumed to be confluent by induction, so there 
exist oriented paths $g_1\to \dots \to g$, $g_2\to \dots \to g$ in 
$\mathcal G_R^{[u]}(X,S)$, and the Diamond Condition holds for $\mathcal G_n^{[v]}$.
\end{proof}

Recall the Shirshov's definition of a composition \cite{Sh62_alg_Lie} in the free Lie algebra 
$\Lie \<U\>$.

Let $f,g\in \Lie \<U\>$ be monic Lie-polynomials, $\bar f = [u]$, $\bar g = [v]$.
We say that $f$ and $g$ form a composition relative to a word $w$ if
$u=u_1u_2$, $v=v_1v_2$, $u_2=v_1$ ($u_i,v_i\in U^*$). Here $w = u_1u_2v_2 = u_1v_1v_2$ 
is a LS-word, and there are two Shirshov braketings:
\[
 \{w_{u\leftarrow *} \} = \{*v_2\}_1, \quad
 \{w_{v\leftarrow *} \} = \{u_1*\}_2.
\]
The Lie polynomial
\[
 (f,g)_w = \{fv_2\}_1 - \{u_1g\}_2
\]
is called a {\em composition of $f$ and $g$ relative to $w$}. 
It is important that 
\begin{equation}\label{eq:CompPrincipal}
 \overline{(f,g)}_w <[w]. 
\end{equation}

It follows from the definition that if $f,g\in S$ then 
\[
[w]\to g_1 = [w] - \{fv_2\}_1
\]
is an edge in $\mathcal G_R(X,S)$, and so is 
\[
[w] \to g_2 = [w] -  \{u_1g\}_2.
\]
Therefore, 
$g_1-g_2 = (f,g)_w$.

Suppose $S$  is a set of monic RLie polynomials such that the rewriting system $\mathcal G_R(X,S)$
is {\em reduced}, i.e., it has the following property:
for every vertex $s\in S$ there is only one edge $s\to 0$ in $\mathcal G_R(X,S)$. We will say 
$S$ is reduced if so is $\mathcal G_R(X,S)$.

\begin{proposition}\label{prop:Ambiguity2}
Let $S$ be a reduced set of monic RLie-polynomials in $\RLie \<X\>$. 
Suppose that all compositions of type $(s_1,s_2)_w$, 
$s_1,s_2\in S$, $[w]\in V_n\cap \RLS(X)$, have the following presentation:
\begin{equation}\label{eq:TrivialComp2}
 (s_1,s_2)_w = \sum\limits_i \alpha_i h_i,\quad h_i\in \hat S\cap V_n,\ \bar h_i<[w].
\end{equation}
Then the rewriting system 
$\mathcal G(V_n)\subset \mathcal G_R(X,S)$ is confluent.
\end{proposition}

\begin{proof}
Check the conditions of Proposition \ref{prop:Ambiguity1} for a word $[w]\in \RLS(X)\cap V_n$.
Assume there is a pair of edges  in $\mathcal G(V_n)$:
  $[w]\to g_1$, $[w]\to g_2$.

There are several possible cases.

1) Both edges are of level 0. 
(This case is actually covered by the classical Composition-Diamond Lemma \cite{Sh62_alg_Lie}, 
but we prefer to consider it in our terminology to make the exposition complete.)
 Then 
 \[
  g_1 = [w] -h_1,\quad g_2 = [w]-h_2,
 \]
$h_i \in \hat \{s_i\}_0\cap V_n$, $s_i\in S$. 
Let $[u]=\bar s_1$ and $[v]=\bar s_2$.

Recall the following

\begin{lemma}[\cite{Chibr}]\label{lem:Nonintersect}
 Suppose $u,v,w\in \LS(U)$, $w=aubvc$, where $a$, $b$, and $c$ are some words in $U$ (either of them may be empty).
 Then there exists a bracketing $\{a*b\star c\}$ such that 
 $\{a[u]b[v] c\} = [w] +\sum\limits_i \alpha_i [w_i]$, $[w_i]<[w]$.
\end{lemma}

This statement also implicitely appears in \cite{Sh62_alg_Lie}.

1.1) Let the corresponding occurences of subwords 
$u$ and $v$ in $w$ do not intersect.
Then $w = aubvc$, $h_1 = \{as_1 bvc\}_1$, $h_2 = \{aub s_2 c\}_2$, 
where $\{\dots\}_1$ and $\{\dots \}_2$ are the Shirshov bracketings on $w_{u\leftarrow *}$ and $w_{v\leftarrow *}$, respectively.
Therefore, 
\[
\begin{aligned}
 g_1-g_2 & = \{aub s_2 c\}_2 - \{as_1 bvc\}_1 \\
 &=\{aub s_2 c\}_2 - \{a[u]bs_2c\}_{12} + \{as_1b[v]c\}_{12} - \{as_1 bvc\}_1 \\
 &+ \{a[u]bs_2c\}_{12} - \{as_1bs_2c\}_{12} + \{as_1bs_2c\}_{12} - \{as_1b[v]c\}_{12} \\
 &= \big(\{aub s_2 c\}_2 - \{a[u]bs_2c\}_{12}\big) + \big(\{as_1b[v]c\}_{12} - \{as_1 bvc\}_1\big ) \\
 &+ \{a([u]-s_1)bs_2c\}_{12} + \{as_1b(s_2-[v])c\}_{12},
\end{aligned}
\]
where $\{a*b\star c\}_{12}$ is the bracketing from Lemma \ref{lem:Nonintersect}.
In the last expression, all summands belong to linear span of $h_i\in \hat S\cap V_n$ with $\bar h_i<[w]$,
so $g_1-g_2$ has the required presentation \eqref{eq:TrivialComp}.

1.2) Let the corresponding occurences of $u$ and $v$ in $w$ intersect:
$u=u_1u_2$, $v=v_1v_2$, $u_2=v_1$. Then $z=uv_2 = u_1v$ is a LS-word,
$w = azb = auv_2b = au_1vb$, 
$h_1 = \{as_1v_2b\}_1$, $h_2 = \{au_1s_2b\}_2$, 
where $\{a*v_2b\}_1$ and $\{au_1*b\}$ are the Shirshov bracketings on $w_{u\leftarrow *}$
and $w_{v\leftarrow *}$, respectively.
Consider also the Shirshov bracketings $\{a*b\}_0$ on $w_{z\leftarrow *}$ and 
$\{*v_2\}_{01}$, $\{u_1*\}_{02}$ on $z_{u\leftarrow *}$ and $z_{v\leftarrow *}$, 
respectively.
Then 
\[
\begin{aligned}
 g_1-g_2 & = \{au_1s_2b\}_2 - \{as_1v_2b\}_1 \\
 &=\{a \{u_1s_2 \}_{02} b\}_0 - \{a\{s_1v_2\}_{01} b\}_0 \\
 &=-\{afb\}_0,
 \end{aligned}
\]
where $f = (s_1,s_2)_z$.
Since $f = \sum\limits_i \alpha_i h_i$, $\bar h_i<[z]$, $h_i\in \hat S\cap V_n$, 
RLie polynomial $g_1-g_2$ may be presented as \eqref{eq:TrivialComp}.

2) The edge 
$[w]\to g_1$ is of positive level $d$, $[w]\to g_2$ is of level 0.

In this case, 
\[
 w = a_1\dots a_m, \quad a_i\in U,
\]
 $a_k = R([v])$ for some $k$, where $[v]=\bar h$, $h\in \widehat{\{s_1\}}_{d-1}$, $s_1\in S$.
Therefore, $h_1 = [a_1\dots a_{k-1}R(h)a_{k+1}\dots a_m]$.
As above,
\[
 w=aub,\quad [u]=\bar s_2,
\]
and $h_2=\{as_2b\}$.
Since $S$ is reduced, the occurence of letter $a_k=R([v])\in U$ considered above 
may appear in either of the subwords $a$ or $b$.
Suppose $a=ca_kc'$ (the second case in analogous).
Then
\[
 w = c R([v])c'ub, \quad c,c',b\in U^*\cup \{\epsilon \} ,
\]
Therefore, 
\[
\begin{aligned}
 g_1-g_2 & = \{cR([v])c's_2b\} - [cR(h)c'ub]     \\
 & = \{cR([v])c's_2b\} -\{cR(h)c's_2b \} +\big (\{cR(h)c's_2b \}- [cR(h)c'ub] \big ),
\end{aligned}
\]
and the same reasonings as in Case 1.1 show the required relation \eqref{eq:TrivialComp} holds.

3) Both edges $[w]\to g_1$, $[w]\to g_2$ have positive level. 
In this case, 
$w = a_1\dots a_k \dots a_l \dots a_m$, $a_i\in U$, where 
$a_k=R([u])$, $a_l = R([v])$, where $[u]\to g_1'$ and $[v]\to g_2'$ are edges of smaller 
level, and 
\[
 h_1  = [a_1\dots R(g_1') \dots a_l \dots a_m], \quad h_2 = [a_1\dots a_k \dots R([v]) \dots a_m].
\]

3.1) If $k\ne l$ then one may proceed as in Case~2.

3.2) If $k=l$, proceed by induction on the level of edges. 
Consider 
$a_k=a_l = R([u])$ with edges $[u]\to g_1'$, $[u]\to g_2'$ in $\mathcal G(V_{n-1})\subset \mathcal G(V_n)$. 
Inductive assumption 
claims  
$g_1'-g_2' = \sum\limits_i \alpha_i h'_i$, $\bar h_i'<[u]$.
Therefore, 
\[
 g_1-g_2 = [a_1\dots R(g_1'-g_2') \dots a_m]
\]
also has a required presentation \eqref{eq:TrivialComp}.
\end{proof}

The entire system $S$ is {\em closed with respect to composition} if for every $s_1,s_2\in S$ every their composition 
$(s_1,s_2)_w$ may be presented as
\[
 (s_1,s_2)_w = \sum\limits_i \alpha _i h_i, \quad h_i\in \hat S,\ \bar h_i <[w],\ \deg_R h_i\le \deg_R w. 
\]
A reduced set of monic RLie polynomials in $\RLie \<X\>$
which is closed with respect to composition is called a {\em Gr\"obner---Shirshov basis} (GSB)
in  $\RLie \<X\>$.

\begin{theorem}\label{thm:CDLemma}
If $S$ is a GSB in $\RLie \<X\>$.
Then the rewriting system $\mathcal G_R(X,S)$ is confluent.
\end{theorem}

\begin{proof}
The statement follows from Propositions \ref{prop:Ambiguity1} and \ref{prop:Ambiguity2}.
\end{proof}

\begin{corollary}
If $S$ is a GSB in $\RLie \<X\>$ then the set of $S$-reduced words 
forms a linear basis of the algebra $\RLie \<X\mid S\>=\RLie \<X\>/I_R(S)$.
\end{corollary}

\begin{proof} 
Terminal vertices of $\mathcal G(X,S)$ are exactly linear combinations of $S$-reduced words.
\end{proof}

\begin{example}
Let $S$ consists of all
$R([u])R([v])$, $[u],[v]\in \RALS(X)$, $[u]>[v]$.
Then $S$ is a reduced system closed with respect to compositions, and 
the set of $S$-reduced words coincides with $\RALS(X)$.
\end{example}

Obviously, $\RLie \< X\mid S \> \simeq \RALie \<X\>$, so $\RALS(X)$ is indeed the linear 
basis of $\RALie \<X\>$.

More general, let $L$ be a Lie algebra, and let $X$ be a linear basis of $L$ which is linearly ordered in some way. 

\begin{example}\label{exm:RA_envelope}
The set $W$ of all words $[w]\in \RALS(X)$ such that $w$ do not contain subwords of type 
$xy$, $x,y\in X$, $x>y$,
form a linear basis of $U_{RA}(L)$.
\end{example}

It is easy to see that $S=\{ R([u])R([w]) \mid [u],[w]\in W, u>w \}\cup \{xy-[x,y]\mid x,y\in X,x>y \}$
is a GSB, and $\RLie \<X\mid S\>\simeq U_{RA}(L)$. 
%
%

\section{Rota---Baxter Lie algebras}

Let $\RBLie $ denotes the variety of Lie algebras equipped with 
a Rota---Baxter operator $R$ of weight $\lambda \in \Bbbk $, 
i.e., a linear map satisfying the following identity:
\[
 [R(x),R(y)] = R([R(x),y]) + R([x,R(y)]) + \lambda R([x,y]).
\]

Consider the forgetful functor $\RBLie \to \Lie $. For every $L\in \Lie $ 
there exists universal enveloping $U_{RB}(L)\in \RBLie $: $L\subset U_{RB}(L)$ 
is a Lie subalgebra, and for every $B\in \RBLie $ and homomorphism 
$\varphi : L\to B$ of Lie algebras there exists unique homomorphism of 
$\RBLie $ algebras $\bar \varphi: U_{RB}(L)\to B$ such that $\bar\varphi|_L = \varphi $.
In this Section, we clarify the structure of $U_{RB}(L)$ and prove an analogue 
of the Poincar\'e---Birkhoff---Witt Theorem.

Suppose $L$ is a Lie algebra with a linear basis $X$. Assume $X$ to be well ordered in some way.
Consider 
\[
S^{(0)} = \{  xy - [x,y] \mid \ x,y\in X,\ x>y \}\subset \Lie \<X\> \subset \RLie \<X\>.
\]
Here $[x,y]$ is a linear form in $X$ equal to the product of $x$ and $y$ in $L$.
Then $S^{(0)}$ is a GSB in $\Lie \<X\>$ and, therefore, in $\RLie \<X\>$. Moreover, 
$L\simeq \Lie \<X \mid S^{(0)} \>$. 

Now, consider 
\[
 \rho(x,y) = R(x)R(y) - R(R(x)y) + R(R(y)x) -\lambda R([x,y]),\quad x,y\in X,\ x>y,
\]
and set $S^{(2)}\subset \RLie\<X\>$ to be the union of $S^{(0)}$ set of all $\rho(x,y)$.
Denote by $\mathcal G^{(2)}$ the subgraph $\mathcal G(V_2)$ of $\mathcal G_R(X,S^{(2)})$. 
Obviously, $S^{(2)}$ is a GSB: it is reduced, and the graph $\mathcal G_R(X,S^{(2)})$
has no ambiguities.

Proceed by induction on $R$-degree. 
Assume a reduced system $S^{(n)}$, $n\ge 2$, is already constructed in such a way that 
the subgraph $\mathcal G^{(n)}=\mathcal G(V_n)\subset \mathcal G_R(X,S^{(n)})$
is a confluent rewriting system. Denote by $T_n$ the set of terminal vertices 
of $\mathcal G^{(n)}$, and let $t_n: V_n\to T_n$ be the linear map that turns 
an $\RLie $ polynomial $f$, $\deg_R f\le n$, into the terminal vertex $t_n(f)$
connected with $f$.

For every two terminal words $a,b \in T_n\cap \RLS(X)$, $\deg_R a + \deg_R b = n-1$, $a>b$, 
consider 
\[
 \rho(a,b) = R(a)R(b) - R(t_n([R(a),b])) + R(t_n([R(b),a])) - \lambda R(t_n([a,b])),
\]
where $[\cdot ,\cdot ]$ stands for the product in $\RLie \<X\>$.
Construct 
\[
 S^{(n+1)} = S^{(n)} \cup \{\rho(a,b) \mid a,b \in T_n\cap \RLS(X),\ \deg_R a + \deg_R b = n-1,\ a>b \}.
\]
It is easy to see from the construction that 
the subgraph $\mathcal G(V_n) \subset \mathcal G_R(X, S^{(n+1)})$
coincides with $\mathcal G^{(n)}$.

We have to resolve the following questions:
\begin{itemize}
 \item Prove that $S^{(n+1)}$ is confluent (assuming so is $S^{(n)}$);
 \item Show  $\RLie \<X\mid S\>\simeq U_{RB}(L)$, where $S$ is the union of all $S^{(n)}$;
 \item Describe the set of $S$-reduced words in $\RLS(X)$.
\end{itemize}

\begin{lemma}\label{lem:Terminal}
 Let $f,g\in \RLie\<X\>$, $\deg_R f + \deg_R g =n$. 
 Then $t_n([f,t_n(g)]) = t_n([f,g])$.
\end{lemma}

\begin{proof}
 It follows from Lemma \ref{lem2} that $[f,g]$ and $[f,t_n(g)]$ belong to the same connected component of $\mathcal G^{(n)}$. 
 Since the latter is confluent, $t_n([f,g]) = t_n([f,t_n(g)])$.
\end{proof}

\begin{lemma}\label{lem:Confluent-n} 
The rewriting system $\mathcal G^{(n+1)}=\mathcal G(V_{n+1})\subset \mathcal G_R(X,S^{(n+1)})$ is confluent. 
\end{lemma}

\begin{proof}
Here we assume by induction that $\mathcal G^{(n)}=\mathcal G(V_n)\subseteq \mathcal G_R(X,S^{(n+1)})$
is confluent.
 It is enough to check the conditions of Proposition \ref{prop:Ambiguity2} 
 for compositions $(s_1,s_2)_w$, $s_1,s_2\in S^{(n+1)}$, $[w]\in \RLS(X)$, $\deg_R w = n+1$.
 
Suppose $s_1 = \rho(a,b)$, $s_2=\rho(b,c)$, $w=R(a)R(b)R(c)$, $a,b,c\in T_n\cap \RLS(X)$, $a>b>c$. 
Denote 
\[
 \begin{gathered}
\rho(a,b) = R(a)R(b) - \sum\limits_i \gamma_i R(c_i), \\  
\rho(b,c) = R(b)R(c) - \sum\limits_j \alpha_j R(a_j), \\
\rho(a,c) = R(a)R(c) - \sum\limits_l \beta_l R(b_l), \\
 \end{gathered}
\]
where 
$\deg_R c_i, \deg_R a_j, \deg_R b_l <n-1$.
Then 
\[
 \begin{aligned}
(s_1,s_2)_w   &  = [\rho(a,b),R(c)]  - [R(a),\rho(b,c)] = [R(a)R(b),R(c)] - [R(a),R(b)R(c)]\\
 &{} - \sum\limits_i \gamma_i [R(c_i), R(c)]  +  \sum\limits_j \alpha_j [R(a),R(a_j)] \\
 & = -[R(b),\rho(a,c)] -  \sum\limits_l \beta_l [R(b),R(b_l)] \\
 &{} - \sum\limits_i \gamma_i [R(c_i), R(c)] +  \sum\limits_j \alpha_j [R(a),R(a_j)] \\
 & = -[R(b),\rho(a,c)] + K(a,b,c).
\end{aligned}
\]
Here $h=[R(b),\rho(a,c)]\in \hat S^{(n)}$, $\hat h = [R(b)R(a)R(c)]<[w]$, and all monomials in 
\[
 K(a,b,c) = \sum\limits_j \alpha_j [R(a),R(a_j)]  -\sum\limits_l \beta_l [R(b),R(b_l)] - \sum\limits_i \gamma_i [R(c_i), R(c)]
\]
are smaller than $[w]$ since they are of degree two in $U$.
Straightforward computations show
\[
 K(a,b,c) = \sum\limits_k \xi_k h_k + R(J(a,b,c)), 
\]
where 
\[
\begin{aligned}\label{eq:R-Jacobian}
 J(a,b,c) 
      &{} = t_n\big ([R(a), t_n([R(b),c] + [b,R(c)] +\lambda [b,c])] + [a, t_n([R(b),R(c)])] \\
         & + \lambda [a, t_n([R(b),c] + [b,R(c)] +\lambda [b,c])] \\
      &{} - [R(b),t_n([R(a),c]+[a,R(c)] +\lambda [a,c])]  - [b,t_n([R(a),R(c)])] \\
         & - \lambda [b, t_n([R(a),c]+[a,R(c)] +\lambda [a,c])] \\
      &{} - [t_n([R(a),R(b)]), c] - [t_n([R(a),b] + [a,R(b)] + \lambda [a,b]), R(c)]\\
         &- \lambda [  t_n([R(a),b] + [a,R(b)] + \lambda [a,b]), c]
      \big ) .
\end{aligned}
\]
Indeed, 
$[R(b),R(c)] \to R(t_n([R(b),c] + [b,R(c)] + \lambda [b,c] )) $
is an edge in $\mathcal G^{(n)}$, so $t_n([R(b),R(c)]) = \sum\limits_j \alpha_j  R(a_j)$.
Moreover, $t_n(R(x)) = R(t_n(x))$ for all $x\in V_{n-1}$.

It remains to apply Lemma~\ref{lem:Terminal} to conclude
\[
\begin{aligned}
 J(a,b,c)  &= t_n\big (
  %
   \Jac(R(a),R(b),c) + \Jac(R(a),b,R(c)) + \Jac(a,R(b),R(c)) \\
   &{} +\lambda \Jac(R(a),b,c) +\lambda \Jac(a,R(b),c) + \lambda \Jac(a,b,R(c)) \\
   &{} +\lambda^2 \Jac(a,b,c)  \big) = 0,
 \end{aligned}
\]
where $\Jac(x,y,z)=[x,[y,z]] - [y,[x,z]] -[[x,y],z]$ is the Jacobian.
Hence, $(s_1,s_2)_w$ has a required presentation \eqref{eq:TrivialComp2}.
\end{proof}

Denote $S = \bigcup\limits_{n\ge 1} S^{(n)}$. Obviously, $S$ is a GSB. Denote by 
$T$ the set of terminal vertices in $\mathcal G_R(X,S)$, 
$T=\bigcup\limits_{n\ge 1} T_n$.
%
%

\begin{lemma}\label{lem:RBaxter}
$\RLie \<X\mid S\>$ is a Rota---Baxter Lie algebra.
\end{lemma}

\begin{proof}
We have to prove 
\begin{equation}\label{eq:RBaxterRelation}
 [R(f),R(g)] - R([R(f),g]) - R([f,R(g)]) - \lambda R([f,g]) \in I_R(S)
\end{equation}
for all $f,g\in \RLie \<X\>$. 
Since for every $f\in \RLie\<X\>$ there exists $t\in T$ such that $f-t\in I_R(S)$,
it is enough to check \eqref{eq:RBaxterRelation} for $f=a$, $g=b$, where 
$a,b\in T\cap \RLS(X)$. 

Assume $a\in T_n$, $b\in T_m$. Then 
$[R(a),R(b)]$ and $R([R(a),b]+[a,R(b)]+\lambda[a,b])$
have the same terminal form in $\mathcal G^{(n+m+2)}$, 
so they are connected by a non-oriented path in $\mathcal G_R(X,S)$. 
Hence, \eqref{eq:RBaxterRelation} holds.
\end{proof}

\begin{corollary}
 $\RLie \<X\mid S\> \simeq U_{RB}(L)$.
\end{corollary}

\begin{proof}
Let $B\in \RBLie$, and let $\varphi : L\to B$ be a homomorphism of Lie algebras.
Identify $L$ with the Lie subalgebra in $\RLie\<X\>$ spanned by $X$. 
Then there exists unique homomorphism of RLie algebras $\psi: \RLie\<X\> \to B$
such that $\psi(x)=\varphi(x)$ for $x\in X$. 
Denote by $\tau $ the natural homomorphism $\RLie\<X\> \to \RLie \<X\mid S\>$, 
$\Ker \tau = I_R(S)$. Since for every $f\in V_n$ we have $f - t_n(f) \in \Ker \tau$, 
Lemma~\ref{lem:RBaxter} implies $S^{(n)}\subset \Ker \tau$. Therefore, 
there exists a homomorphism of RLie algebras $\bar\varphi:\RLie\<X\mid S\> \to B $, 
$\bar\varphi(x)=\psi(x) = \varphi(x)$ for $x\in X$.
\end{proof}

The universal enveloping Rota---Baxter Lie algebra $U_{RB}(L)$ of a Lie algebra $L$ 
has a natural ascending filtration induced by
\eqref{eq:Deg-Filtration}:
\[
 U_{RB}^{(n)}(L) = \tau(\RLie^{(n)}\<X\>).
\]
Denote by $\gr U_{RB}(L) $ the associated graded RLie algebra. The following statement is ideologically similar
to the classical Poincar\'e---Birkhoff---Witt Theorem.

\begin{theorem}\label{thm:PBW_Lie}
$\gr U_{RB}(L) \simeq U_{RA}(L)$ as  Lie algebras. 
\end{theorem}

\begin{proof}
It is enough to compare Gr\"obner---Shirshov bases of $U_{RA}(L)$ 
and 
$U_{RB}(L)$. The principal parts of these relations coincide, they are of degree~2. For the latter algebra, 
the right-hand sides of relations are of degree~1. 
\end{proof}

\begin{remark}
For associative algebras the statement of Theorem \ref{thm:PBW_Lie} is easy 
to show by means of Gr\"obner---Shirshov bases technique for associative Rota---Baxter algebras \cite{Bokut2010}: 
multiplication table of an associative algebra $A$ is 
closed under all compositions in the free associative Rota---Baxter algebra.
\end{remark}

\begin{remark}
 For commutative algebras, an analogue of Theorem \ref{thm:PBW_Lie} also holds. Moreover, 
 there is an explicit construction of the universal 
 enveloping commutative Rota---Baxter algebra $U_{RB}(A)$ for a given commutative algebra~$A$ 
 (mixed shuffle algebra \cite{Guo_RBA}). 
\end{remark} 
 
Let us briefly state the construction from \cite{Guo_RBA} (in the case of zero weight) in more natural terms.
Consider 
 \[
 \mbox{\CYRSH}\, (A) = A^\# \otimes B^\#, \quad B=\preCom\<A^\#\>^{(+)}, 
 \]
where 
$\preCom\<A^\#\>$ is the free pre-commutative (Zinbiel) algebra generated by the space $A^\#$, 
$B$ is the anti-commutator algebra of $Z$ (it is an associative and commutative algebra), and
$A^\#=A\oplus \Bbbk 1_A$,
$B^\#=B\oplus \Bbbk 1_B$
are obtained by joining external units.

Define the linear operator on $\mbox{\CYRSH}\, (A)$:
\[
 R(a\otimes 1_B) = 1_A\otimes a, \quad a\in A^\#,
\]
\[
 R(a\otimes b ) = 1_A\otimes ab,\quad a\in A^\#, \ b\in B.
\]
The Zinbiel identity $(xy)z=x(yz)+x(zy)$ on $\preCom \< A^\# \>$
implies $R$ to be a Rota---Baxter operator on $\mbox{\CYRSH}\, (A)$. For example,
\begin{multline}\nonumber
 R(a_1\otimes 1_B)R(a_2\otimes b) = (1_A\otimes a_1)(1_A\otimes a_2b) \\
 = 1_A\otimes (a_1(a_2b) + (a_2b)a_1) = 1_A\otimes a_1(a_2b) + 1_A\otimes a_2(ba_1) + 1_A\otimes a_2(a_1b).
\end{multline}
On the other hand,
\begin{multline}\nonumber
 R((a_1\otimes 1_B)R(a_2\otimes b) +R(a_1\otimes 1_B)(a_2\otimes b ) )
 =
 R( (a_1\otimes 1_B)(1_A \otimes a_2b) +(1_A\otimes a_1)(a_2\otimes b )  )\\
 =R( a_1\otimes a_2b  + a_2\otimes( a_1b +ba_1 )  )
 =1_A\otimes a_1(a_2b) + 1_A\otimes a_2(a_1b) + 1_A\otimes a_2(ba_1).
\end{multline}

It is easy to check (see \cite{Guo_RBA}) that the embedding 
\[
 A\to \mbox{\CYRSH}\, (A),\quad a\mapsto a\otimes 1_B, \ a\in A,
\]
may be extended to a homomorphism of Rota---Baxter algebras
$U_{RB}(A)\to \mbox{\CYRSH}\, (A)$. 
Suppose $X$ is a linear basis of $A$ and consider the  following elements of $U_{RB}(A)$: 
\begin{gather}
u = R^{s_1}(x_1R^{s_2}(x_2R^{s_3}(\dots R^{s_{n-1}}(x_{n-1}R^{s_n}(x_{n}))\dots )),             \label{eq:ShuffleBasis}
\\  
x_i\in X,\ s_1\ge 0, \ s_2,\dots s_n>0, n\ge 1.\nonumber
\end{gather}
Images of these elements $\mbox{\CYRSH}\, (A)$ are linearly independent since
the linear base of $\preCom\<A^\#\>$ is given by $x_1(x_2(\dots (x_{n-1}x_n)\dots ))$, $x_i\in X\cup\{1_A\}$ \cite{Loday01}.
On the other hand, the set of \eqref{eq:ShuffleBasis} obviously span $U_{RB}(A)$.
Therefore, \eqref{eq:ShuffleBasis} is a linear basis of $U_{RB}(A)$ as well as of $U_{RA}(A)$ from Remark~\ref{rem:RAVar}.

For nonzero weight, it is enough to 
replace $\preCom\<A^\#\>$ with $\mathop{\fam 0 postCom}\,\<A^\#\>$ (commutative tridendriform algebra, or CTD-algebra), 
and set $B$ to be the associated commutative algebra \cite[p.~26]{Zinbiel}.

\begin{remark}
The same statement holds for algebras with a Nijenhuis operator, i.e., a linear map 
$N$ such that 
\[
 [N(x),N(y)] = N([N(x),y]) + N([x,N(y)]) - N^2([x,y]).
\]
\end{remark}

The route of the proof is completely similar to stated above.
The key computation of a composition is based on the following relation which is easy to check 
by straightforward computation:
\[
\begin{aligned}
 \Jac(N(a),N(b),& N(c)) 
  = N\big (\Jac(a,N(b),N(c))+\Jac(N(a),b,N(c)) \\
 &{} +\Jac(N(a),N(b),c)\big) - N^2\big ( \Jac(a,b,N(c)) +\Jac(a,N(b),c) \\
 &{} +\Jac(N(a),b,c)\big ) + N^3\big (\Jac(a,b,c)\big ).
\end{aligned}
\]

\end{document}